\documentclass[10pt,reqno]{article}

\usepackage[b5paper,top=1.0in,left=0.9in]{geometry}
\usepackage{verbatim}
\usepackage{amsthm}
\usepackage{amssymb}
\usepackage{amsmath}
\usepackage{graphicx}
\usepackage{appendix}
\usepackage{color} 
 \usepackage{float}

\newcommand{\D}{\mathrm{D}}

\newtheorem{Thm}{Theorem}[section]
\newtheorem{Lem}[Thm]{Lemma}
\newtheorem{Prop}[Thm]{Proposition}

\newtheorem{Def}[Thm]{Definition}

\newtheorem{Fact}[Thm]{Fact}
\newtheorem{Nota}[Thm]{Notation}

\def\R{\mathbb{R}}

\def\C{\mathbb{C}}
\def\Z{\mathbb{Z}}

\def\to{\longrightarrow}

\def\cM{\mathcal{M}}

\def\cU{\mathcal{U}}

\def\G{\Gamma}
\def\c{\gamma}
\def\D{\Delta}

\def\h{\theta}

\def\t{\tau}
\def\W{\Omega}
\def\w{\omega}
\def\ze{\zeta}

\def\sl{\mathfrak{sl}}

\def\o+{\oplus}
\def\bo+{\bigoplus}
\def\x{\times}
\def\p[#1,#2]{\phi_{#1,#2}}
\def\til[#1]{\widetilde{#1}}
\def\what[#1]{\widehat{#1}}

\def\z[#1]{z_{#1}}

\def\oo{\infty}
\def\=>{\Longrightarrow}

\def\<{\langle}
\def\>{\rangle}

\def\^{\wedge}
\def\+{\dagger}

\def\inv{^{-1}}
\def\dis{\displaystyle}
\def\over[#1]{\overline{#1}}
\def\vec[#1]{\overrightarrow{#1}}
\def\mat[#1, #2]{\left[\begin{array}{ccccc}#1\end{array}\left|\begin{array}{c}#2\end{array}\right.\right]}
\def\xto[#1]{\xrightarrow{#1}}
\def\dd[#1,#2]{\frac{d#1}{d#2}}
\def\del[#1,#2]{\frac{\partial #1}{\partial #2}}
\def\Facts[#1]{\begin{Fact}\mbox{}\begin{itemize}#1\end{itemize}\end{Fact}}
\def\Notation[#1]{\begin{Nota}\mbox{}\begin{itemize}#1\end{itemize}\end{Nota}}
\def\Eqn[#1]{\begin{eqnarray*}#1\end{eqnarray*}}
\def\tab{\;\;\;\;\;\;}

\newcommand{\Eq}[1]{\begin{align}#1\end{align}}

\numberwithin{equation}{section}

\begin{document}
\title{The Graphs of Quantum Dilogarithm}

\author{  Ivan C.H. Ip\footnote{
          Yale University,
          Department of Mathematics,
          10 Hillhouse Avenue,
          New Haven, CT 06520,
          U.S.A.\newline
          Email: ivan.ip@yale.edu}
          }

\date{\today}

\numberwithin{equation}{section}

\maketitle

\begin{abstract}
Using the complex coloring method, we present the graphs of the quantum dilogarithm function $G_b(z)$ and visualize its analytic and asymptotic behaviors. In particular we demonstrate the limiting process when the modified $G_b(z)\to \G(z)$ as $b\to 0$. We also survey the relations of $G_b(z)$ with different variants of the quantum dilogarithm function.
\end{abstract}


\tableofcontents
\section{Introduction} Let $q=e^{\pi ib^2}$ be the quantum parameter. The quantum dilogarithm function $G_b(z)$, defined for $b\in \R_{>0}$, is also known as the quantum exponential function or the hyperbolic gamma function. It bears its name from the properties of its variants $g_b(z)$ and $\Phi_b(z)$ given by
\Eq{\Phi_b(z)=g_b(e^{2\pi bz}) = \frac{\over[\ze_b]}{G_b(\frac{Q}{2}-iz)},}
where $Q=b+b\inv$ and $\ze_b=e^{\frac{\pi i}{2}(\frac{b^2+b^{-2}}{6}+\frac{1}{2})}$,
such that
\Eq{\label{qexp}g_b(u)g_b(v)=g_b(u+v),}
\Eq{\label{qpenta}g_b(v)g_b(u)=g_b(u)g_b(q\inv uv)g_b(v)}
hold for the positive self-adjoint Weyl type operators $u,v$ with $uv=q^2vu$ in the case $|q|=1$. The first relation says that it is a $q$-deformation of the exponential function, while the second relation is the deformed version of Roger's pentagon identity for the classical dilogarithm function.
Using $\Phi_b(z)$, the pentagon relation can be written formally for the Heisenberg type operator $[p,x]=~\frac{1}{2\pi i}$ as well by
\Eq{\Phi_b(p)\Phi_b(x)=\Phi_b(x)\Phi_b(p+x)\Phi_b(p),}
where we ignored the concerns of unboundedness and self-adjointness of operators.

On the other hand, the appearance of $G_b(z)$ in the $q$-binomial formula \cite{BT, Ip1} given by
\Eq{(u+v)^{it} = b\int_{\R+i0}\frac{G_b(ib\t-ibt)G_b(-ib\t)}{G_b(-ibt)}u^{i \t}v^{i (t-\t)} d\t,}
where $uv=q^2vu$ as above, suggests that $G_b(bz)$ serves as the quantum analogue of the classical Gamma function. In fact it is shown in \cite{Ip1} that we have the limit given by
\Eq{\label{Glim}\lim_{r\to 0} \frac{G_b(bz)}{\sqrt{-i}|b|  (1-q^2)^{x-1}} = \G(z),}
where we analytic continued $b$ so that $b^2=ir$ with $r>0$ and the principal values of the complex powers are taken.

Faddeev's quantum dilogarithm function $\psi(z)$ is originally defined for $0<q<1$ using an infinite product expansion \cite{Fa}, and it is shown that it carries an analytic continuation to the regime $|q|=1$, where the infinite product no longer converges, but instead an integral expression is found \cite{FKa}. Since then a lot of variants have been defined in the literature serving different purposes. We have for example the hyperbolic gamma function $\c(z)$ by Volkov \cite{Vo}, the quantum dilogarithm function $\Phi^\h(z)$ by Fock-Goncharov \cite{FG}, $V_\h(z)$ by Woronowicz \cite{WZ}, and the 2-parameter G-function $G(a_+,a_-;z)$ given by Ruijinaars \cite{Ru}. This paper atempts to relate all these functions to the present choice of Teschner's definition of $G_b(z)$ in \cite{BT} so that the graphical presentations for this function can be translated to all these variants by shifting, stretching and phase changes easily.

The quantum dilogarithm function played a prominent role in the representation theory of noncompact quantum groups. This function and its many variants are being studied \cite{Gon, KaN, Ru, Vo} and applied to vast amount of different areas, for example the construction of the '$ax+b$' quantum group by Woronowicz et.al. \cite{PuW, WZ}, the harmonic analysis of the non-compact quantum group $U_q(\sl(2,\R))$ and its modular double \cite{BT,PT1,PT2}, the $q$-deformed Toda chains \cite{KLS} and hyperbolic knot invariants \cite{Ka2}. Recently attempts have also been made to cluster algebra \cite{FG, KaN} and quantization of the Teichm\"{u}ller space \cite{CF, FK, Ka3}. One of the important properties of this function is its invariance under the duality $b\leftrightarrow b\inv$ that help encodes the detail of the modular double of the quantum plane \cite{Ip2}, and also relates, for example, to the self-duality of Liouville theory \cite{PT1}. 

The present paper is organized as follows. In Section \ref{sec:complex}, we introduce Jan Homann's representation of a complex valued function on the complex plane by a coloring method using a Mathematica code, and presents several examples illustrating this idea. Then in Section \ref{sec:qdilog} we define $G_b(z)$ and describe its main analytic and asymptotic properties. In Section \ref{sec:graphs} we present the complex graphs $G_b(z)$ for several values of $b$, and its symmetric version $S_b(z)$. Then we restrict to the most interesting case along the real line and the imaginary direction and study its analytic and asymptotic properties. In Section \ref{sec:compact}, we present the graphs for the compact quantum dilogarithm corresponding to $Re(b^2)=0$, and demonstrate using the complex graphs the limiting process of the modified $\til[G_b](z)$ tending towards the classical gamma function $\G(z)$. Finally in Section \ref{sec:relation}, we relate $G_b(z)$ to all other variants of the quantum dilogarithm function used in the literature. 

\textbf{Acknowledgements.} I would like to thank Alvin Wong and Hyun Kyu Kim for helpful discussions.

\section{Complex Coloring}\label{sec:complex}
\subsection{Specifications} We adopt the coloring method for a complex function developed by Jan Homann. The color function is obtained from the Mathematica code by Axel Boldt and can be found in \cite{AB}. 

In short, given a complex valued function defined on the complex plane, we represent the function by coloring the complex plane according to the following rules. The argument of a complex value is encoded by the hue of a color (red = positive real, and then counterclockwise through yellow, green, cyan, blue and purple; cyan stands for negative real). The absolute value of the function is represented by the brightness and saturation of the color: strong colors denote points close to the origin, with black = 0, weak colors denote points with large absolute value, with white = $\oo$.

The color is defined by 3 numbers ranging from 0 to 1. In the present coding, at the point $z=x+iy$ with value $f=f(x+iy)$, they are explicitly given by
\Eq{Hue[h,s,b]:=\{\frac{Arg(f)}{2\pi}, \frac{1}{1+0.3 \ln(|f|+1)}, 1-\frac{1}{1.1+5\ln(|f|+1)}\},}
where $h$ stands for hue, $s$ the saturation and $b$ the brightness.

\begin{figure}[H]
\centering
\includegraphics[width=60mm]{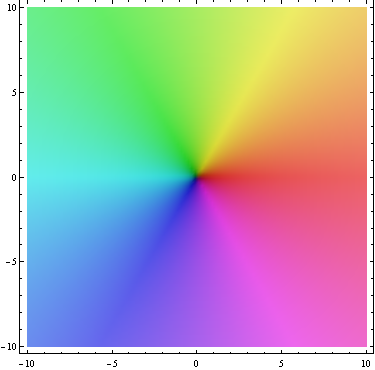}
\includegraphics[width=60mm]{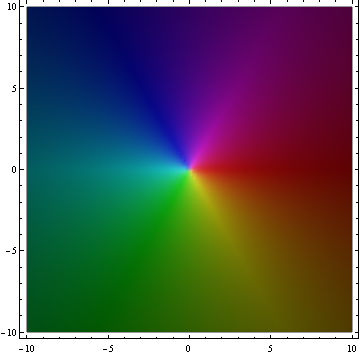}
\caption{The color function for $f(z)=z$ and $f(z)=1/z$}
\label{color}
\end{figure}

The color function is hence given by Figure \ref{color}. We also note that taking inverse amounts to flipping the brightness of the color, and taking the opposite hue along the imaginary direction.

\begin{figure}[htp]
\centering
\includegraphics[width=60mm]{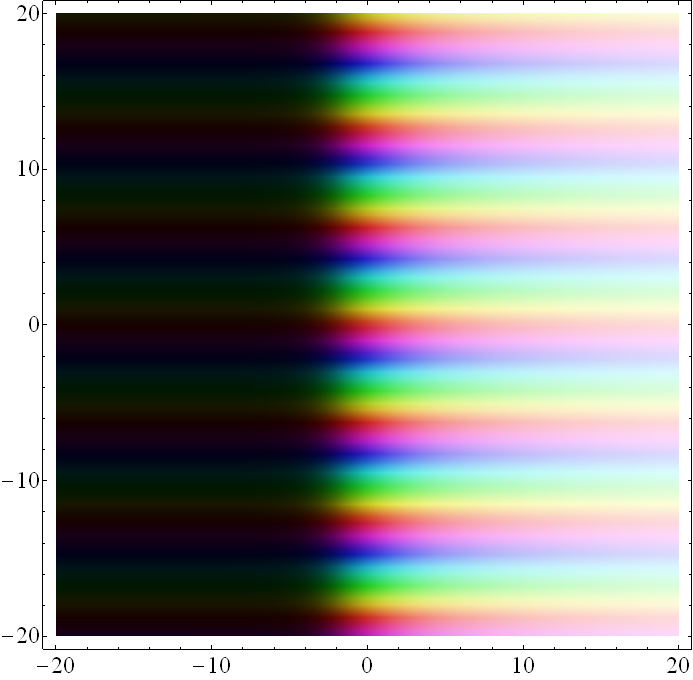}
\caption{The brightness function for $f(z)=e^z$}
\label{exp}
\end{figure}

A good reference for the brightness of the color is given by the exponential function $f(z)=e^z$, see Figure \ref{exp}.

\subsection{Examples}
Let us demonstrate the use of the complex coloring method. We illustrate here the graph of $f(z)=z^3$, the gamma function $f(z)=\G(z)$ and the Riemann Zeta function $f(z)=\ze(z)$.
\begin{figure}[H]
\centering
\includegraphics[width=60mm]{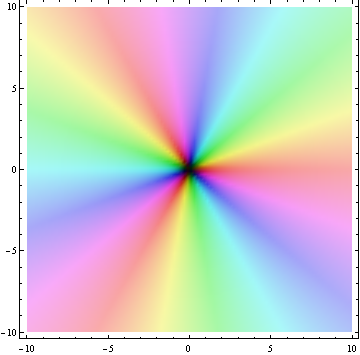}
\caption{The graph of $f(z)=z^3$}
\label{z3}
\end{figure}

From Figure \ref{z3}, we observe that for $f(z)=z^3$ it has an order 3 zero, with the hue of color cycling around the point $z=0$ for 3 times.

\begin{figure}[H]
\centering
\includegraphics[width=50mm]{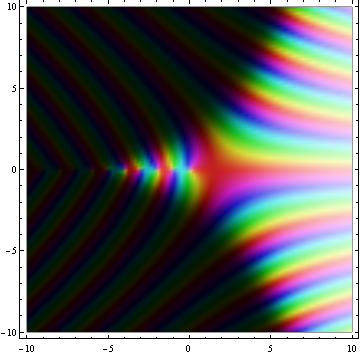}
\caption{The graph of $f(z)=\G(z)$}
\label{Gamma}
\end{figure}

From Figure \ref{Gamma} for $f(z)=\G(z)$, we observe the phase change behavior and the exponential growth along the imaginary direction. The few bright spot at $z=-n,n\in\Z_{\geq0}$ indicates the simple poles.

\begin{figure}[H]
\centering
\includegraphics[width=50mm]{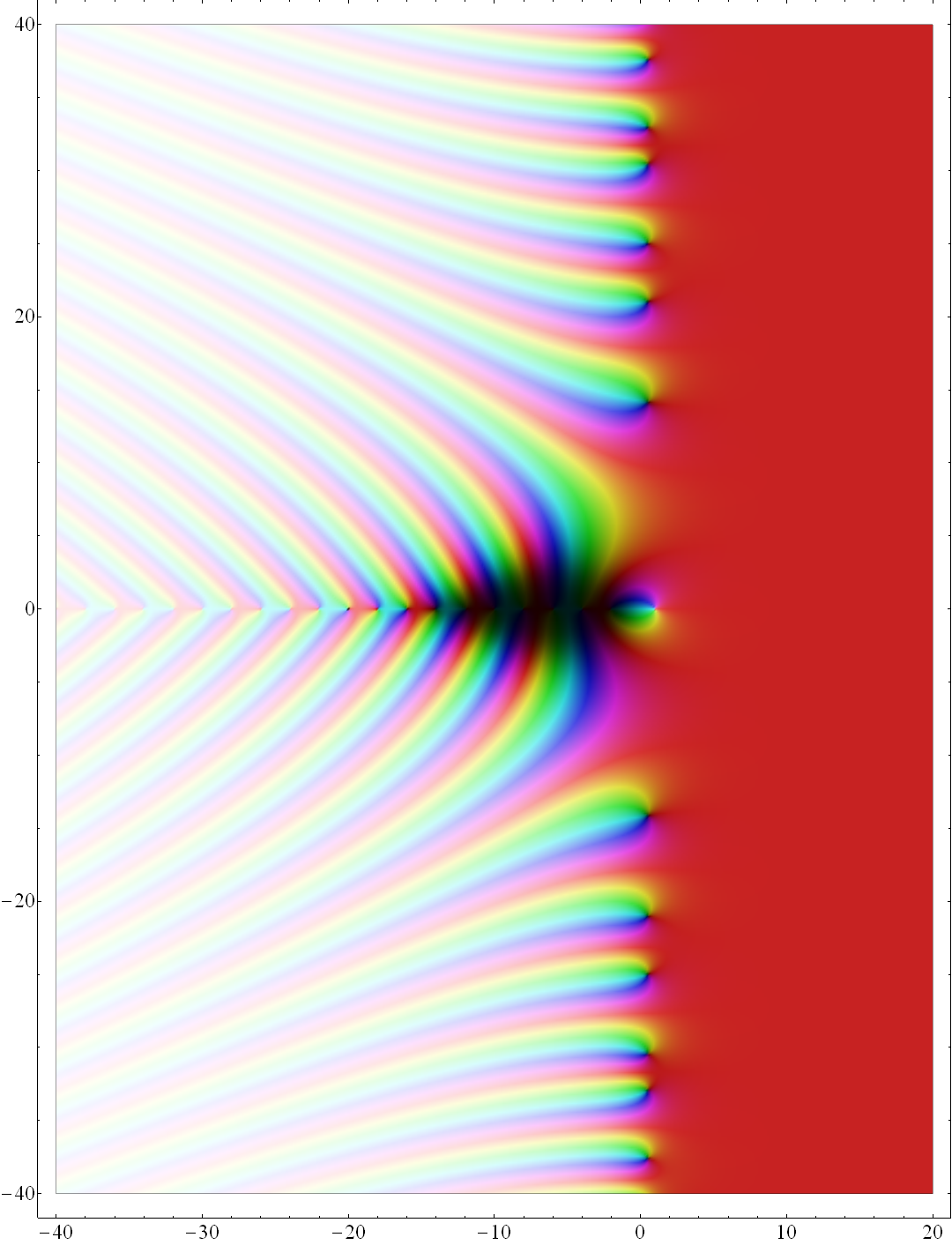}
\caption{The graph of $f(z)=\ze(z)$}
\label{zeta}
\end{figure}

From Figure \ref{zeta} for $f(z)=\zeta(z)$, we observe the existence of nontrivial zeros along the line $Re(z)=~\frac{1}{2}$, the trivial zeroes at $z=-2n,n\in\Z_{\geq1}$ (where the color lines come together), the pole at $z=1$, and the limit $\dis\lim_{Re(z)\to+\oo}\ze(z)=1$.

\section{The Quantum Dilogarithm $G_b(z)$} \label{sec:qdilog}

The quantum dilogarithm function can be defined in two ways. Historically it came from Barnes' double Zeta function \cite{Ba} which is defined as follows:
\begin{Def}Let $\w:=(w_1,w_2)\in \C^2$. The double Zeta function is defined as
\Eq{\ze_2(s,z|\w):=\sum_{m_1,m_2\in\Z_{\geq0}}(z+m_1w_1+m_2w_2)^{-s}.}

The double Gamma function is defined as
\Eq{\G_2(z|\w):=\exp\left(\frac{\partial}{\partial
s}\ze_2(s,z|\w)|_{s=0}\right).}

Let
\Eq{\G_b(x):=\G_2(x|b,b\inv).}

The quantum dilogarithm is defined as the function:
\Eq{S_b(x):=\frac{\G_b(x)}{\G_b(Q-x)},}
where $Q=b+b\inv$.

The main object of study is the variant given by
\Eq{\label{GS}G_b(x):=e^{\frac{\pi i}{2}x(x-Q)}S_b(x).}
\end{Def}

However for computational purpose, it is easier to use the equivalent integral form
\begin{Prop} For $0\leq Re(z)\leq Q$, $G_b(z)$ can be expressed as
\Eq{\label{intform} G_b(x)=\over[\ze_b]\exp\left(-\int_{\W}\frac{e^{\pi tz}}{(e^{\pi bt}-1)(e^{\pi b\inv t}-1)}\frac{dt}{t}\right),}
where \Eq{\ze_b=e^{\frac{\pi i}{2}(\frac{b^2+b^{-2}}{6}+\frac{1}{2})},}
and the contour goes along $\R$ with a small semicircle going above the pole at $t=0$.

The integral converges absolutely for $0<Re(z)<Q$ and conditionally at the boundary.
\end{Prop}

The function $G_b$ and $S_b$ can be analytic continued to the whole complex plane as a meromorphic function, so that they have simple zeros at $z=Q+nb+mb\inv$ and simple poles at $z=-nb+-mb\inv$ for $n,m\in \Z_{\geq 0}$.

The quantum dilogarithm satisfies the following properties:
\begin{Prop} Self-duality:
\Eq{S_b(x)=S_{b\inv}(x),\tab G_b(x)=G_{b\inv}(x).}

Functional equations: \Eq{\label{funceq}S_b(x+b^{\pm1})=2\sin(\pi
b^{\pm1}x)S_b(x),\tab G_b(x+b^{\pm 1})=(1-e^{2\pi ib^{\pm 1}x})G_b(x).}

Reflection property:
\Eq{\label{reflection}S_b(x)S_b(Q-x)=1,\tab G_b(x)G_b(Q-x)=e^{\pi
ix(x-Q)}.}

Complex conjugation: \Eq{\overline{S_b(x)}=\frac{1}{S_b(Q-\over[x])},\tab \overline{G_b(x)}=\frac{1}{G_b(Q-\bar{x})},}

in particular 
\Eq{\left|S_b(\frac{Q}{2}+ix)\right|=\left|G_b(\frac{Q}{2}+ix)\right|=1 \mbox{ for $x\in\R$}.}

\label{asymp} Asymptotic properties:
\Eq{G_b(x)\sim\left\{\begin{array}{cc}\bar{\ze_b}&Im(x)\to+\oo,\\\ze_b
e^{\pi ix(x-Q)}&Im(x)\to-\oo,\end{array}\right.}
where
\Eq{\ze_b=e^{\frac{\pi i}{4}+\frac{\pi i}{12}(b^2+b^{-2})}.}

\label{residue} Residues:
\Eq{\lim_{x\to 0} xG_b(x)=\frac{1}{2\pi},} or more generally,
\Eq{Res\frac{1}{G_b(Q+z)}=-\frac{1}{2\pi}\prod_{k=1}^n(1-q^{2k})\inv\prod_{l=1}^m(1-\widetilde{q}^{2l})\inv}
at $z=nb+mb\inv, n,m\in\Z_{\geq0}$ and $\widetilde{q}=e^{\pi i
b^{-2}}$.
\end{Prop}

Using the functional equation \eqref{funceq}, the function can be extended to the whole complex plane using the integral formula \eqref{intform} by
\Eq{G_b(z+nQ)=G_b(x)\prod_{k=1}^n (1-e^{2\pi ib(z+(k-1)Q)}),}
\Eq{G_b(z-nQ)=G_b(x)\prod_{k=1}^n (1-e^{2\pi ib(z-kQ)})\inv,}
where $n\in \Z_{\geq 0}$, $0\leq Re(z)\leq Q$.

\section{Visualizing $G_b(z)$ and $S_b(z)$}\label{sec:graphs}
\subsection{On the complex plane}
Using the complex coloring method developed in Section \ref{sec:complex}, we obtained the graph for $G_b(x)$ for 2 random numbers $b=0.7$ and $b=0.248$ which suffice in showing the general shape of the function. We also considered the limiting case for $b=1$. Due to slow convergence in the numerical calculation, we will consider another limiting case as $b=0.1$ only in the next subsection when we restrict to the real line.

\begin{figure}[H]
\centering
\includegraphics[width=140mm]{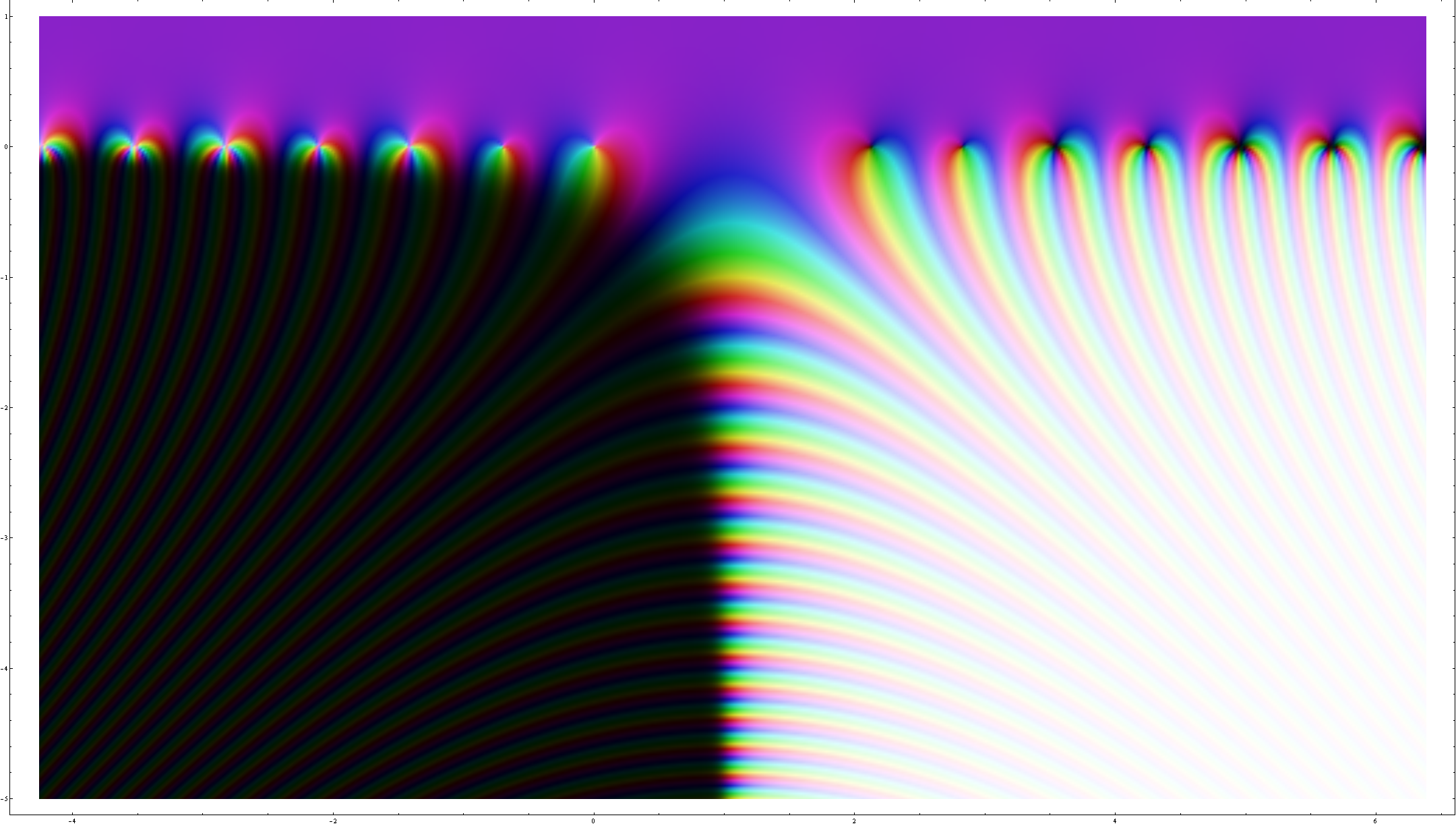}
\caption{The graph of $G_b(z)$ for $b=0.7, z\in [-2Q,3Q]\x[-5,1], Q=2.129$}
\label{Gx07}
\end{figure}

\begin{figure}[H]
\centering
\includegraphics[width=140mm]{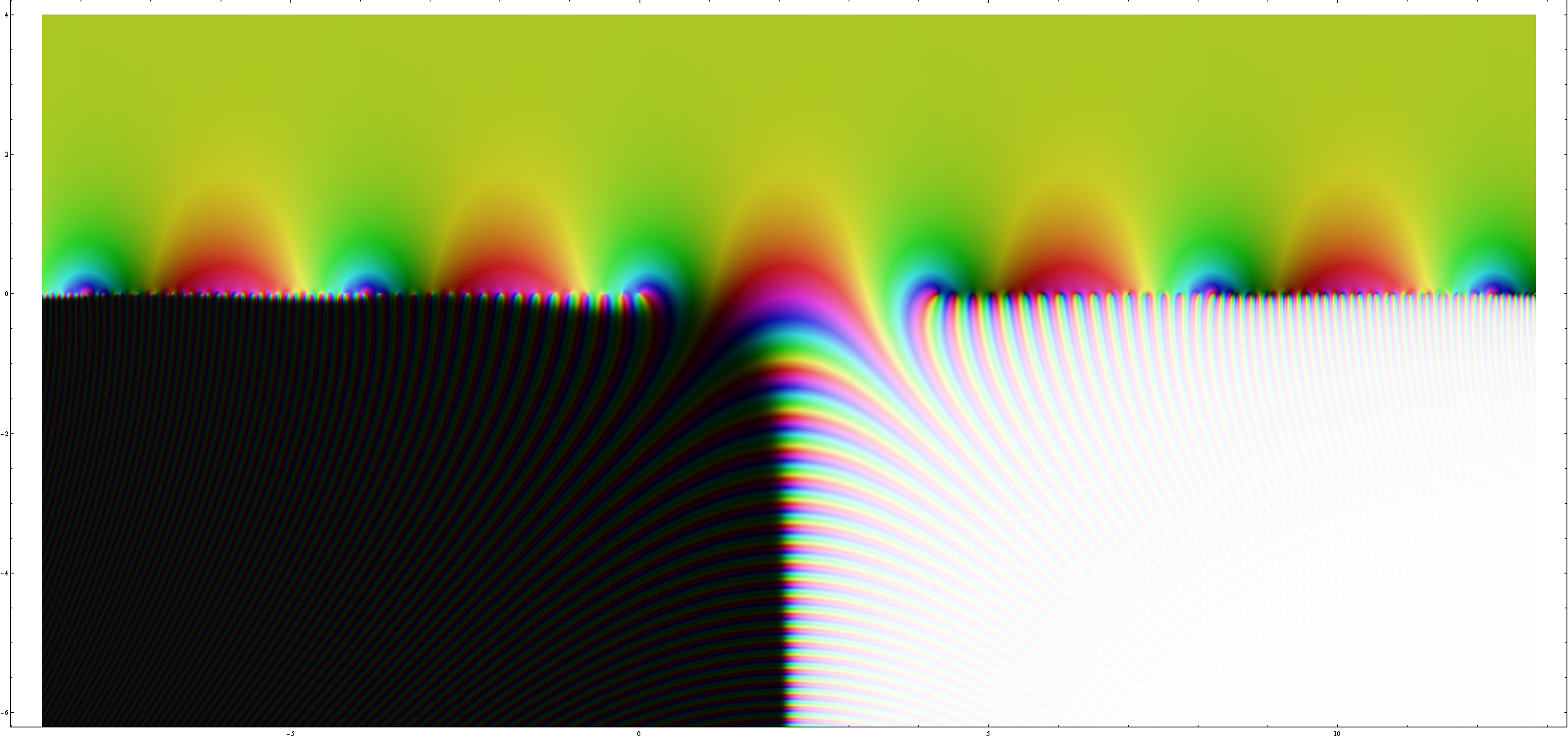}
\caption{The graph of $G_b(z)$ for $b=0.248, z\in [-2Q,3Q]\x[-6,4], Q=4.280$}
\label{Gx0248}
\end{figure}

The graph of $G_b(z)$ for $b=0.7$ is given in Figure \ref{Gx07}, while the graph for $b=0.248$ is given in Figure \ref{Gx0248}.

We immediately observe several properties
\begin{itemize}
\item We observe the symmetry along the axis $Re(z)=\frac{Q}{2}$, governed by the reflection properties \eqref{reflection}. Furthermore, the poles and zeros along $Im(z)=0$ are supposed to be simple, by comparing with the graph of $f(z)=z$ given in Section \ref{sec:complex}. The appearance of certain double or triple poles and zeros for the graph $b=0.7$ comes from the fact that $b\inv\simeq1.43\sim 2b$, hence certain multiples $nb+mb\inv$ will come close together. For example, $z=2b+2b\inv=4.257$ and $z=4b+b\inv=4.229$.
\item It has exponential decay along the negative imaginary direction for $Re(z)<\frac{Q}{2}$, and exponential growth for $Re(z)>\frac{Q}{2}$. According to the asymptotic behavior from the last section, the growth is given by 
\Eq{|G_b(s+it)| = e^{-\pi t(2s-Q)},\tab t\to -\oo.}
This also explains the brightness for different $b$ when $Im(z)\to -\oo$.
\item On the other hand, it approaches the limit $\overline{\ze}_b = e^{\pi i(b^2+b^-2+3)/12}$ in the positive imaginary direction.
\item We observe the phase change along any given direction. In particular along $Re(z)=\frac{Q}{2}$ we do have the changes given according to $e^{\pi iz^2}$, increasing in frequency as $Im(z)\to -\oo$.
\item Finally let us note the similar shape compared with the Riemann Zeta function, which we believe is a coincidence.
\end{itemize}

It is also interesting to look at the limiting case for $b=1$.

\begin{figure}[H]
\centering
\includegraphics[width=140mm]{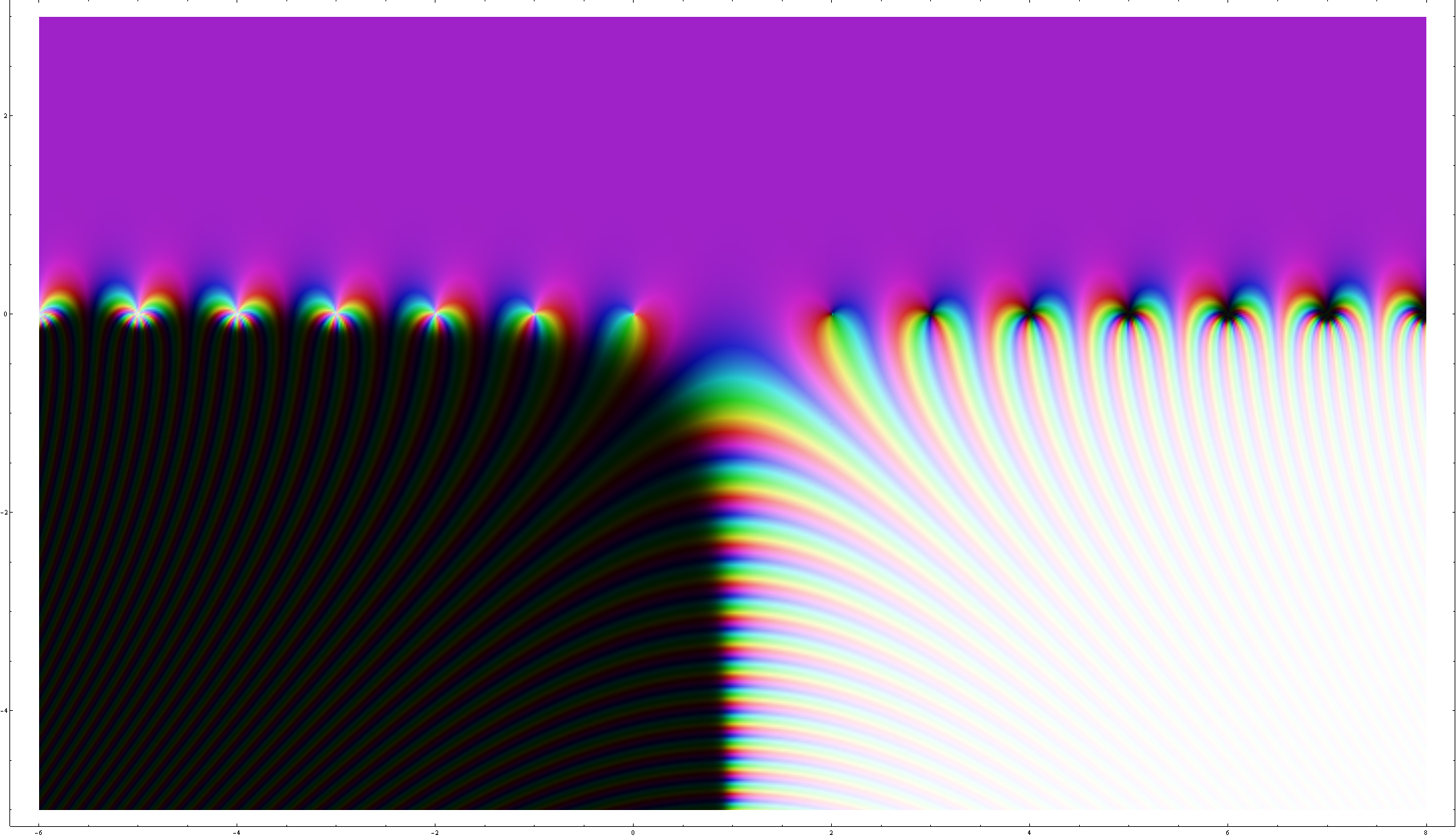}
\caption{The graph of $G_b(z)$ for $b=1, z\in [-2Q,3Q]\x[-5,3], Q=2$}
\label{Gx1}
\end{figure}

For the case $b=1$, from Figure \ref{Gx1} we observe that the successive poles and zeroes have increasing order, due to the fact that $b$ and $b\inv$ coincide, which creates multiple poles. The poles at $z=-n$ and the zeroes at $z=2+n$ both have order $n+1$ for $n=0,1,2,...$

\begin{figure}[H]
\centering
\includegraphics[width=140mm]{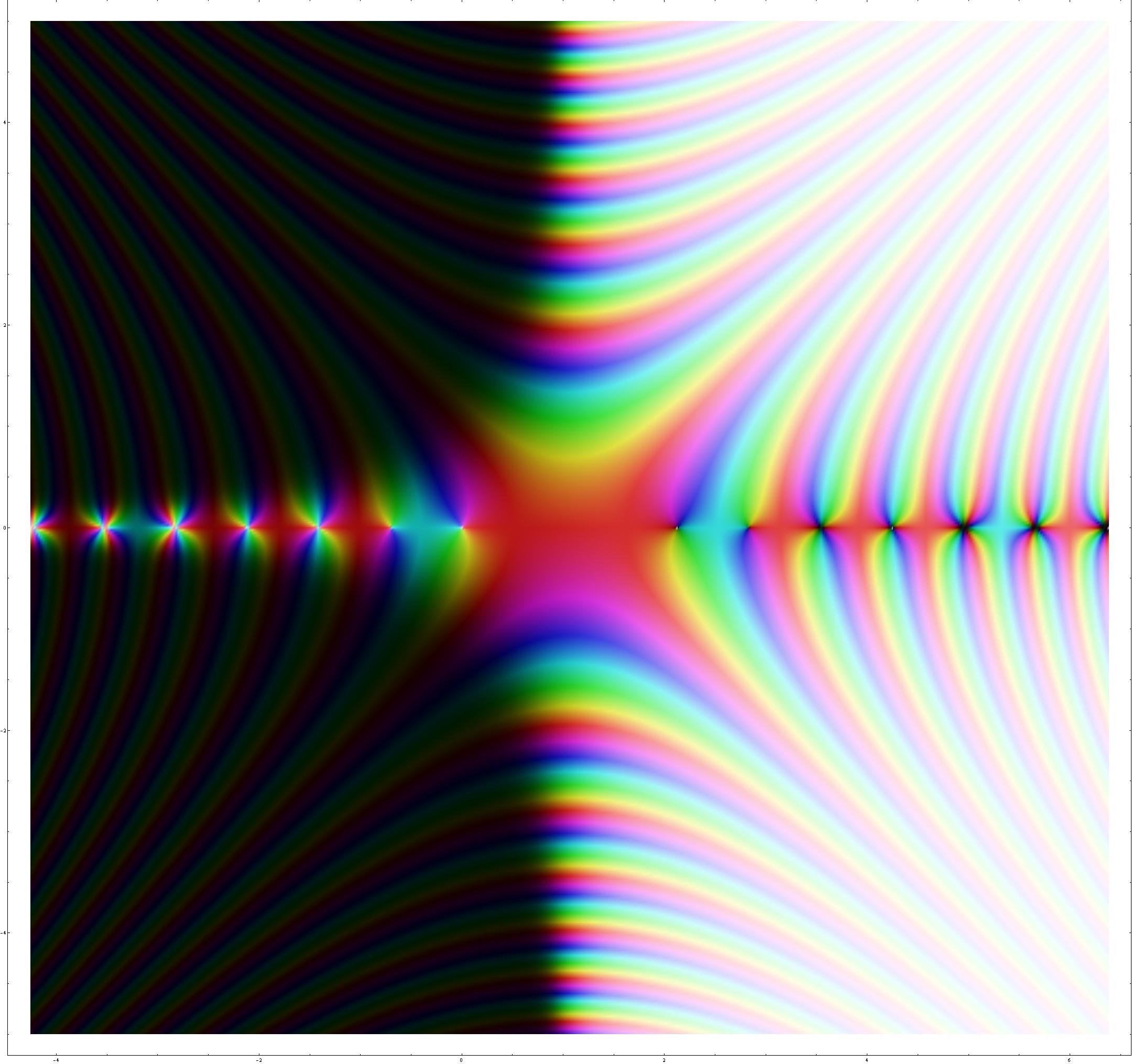}
\caption{The graph of $S_b(z)$ for $b=0.7, z\in [-2Q,3Q]\x[-5,5], Q=2.129$}
\label{Sx}
\end{figure}

The graph can be made symmetric by considering $S_b(z)$ instead. The reflection properties give
\Eq{S_b(z)=S_b(Q-z)\inv,\tab \over[S_b(z)]=S_b(\over[z]),}
suggesting a full symmetry of the graph along $Re(z)=\frac{Q}{2}$ and $Im(z)=0$. We illustrate the case for $b=0.7$ in Figure \ref{Sx}.

\subsection{Values along $z\in\R$}
The most interesting behavior of the quantum dilogarithm appears at the real line, where all the poles and zeroes lie. Since we know that $\over[S_b(z)]=S_b(\over[z])$, we conclude that $S_b(z)$ takes real values when $z$ is real. Therefore it suffices to plot the graphs for $S_b(z)$, and the values for $G_b(z)$ just amount to a phase change by the factor $e^{\pi iz(z-Q)/2}$ (cf. \eqref{GS}).

Here we plot the graph for $b=0.7, b=0.248$ as well as the limiting case $b=0.1$. Mathematica obtains a numerical overflow beyond this small value of $b$ due to large exponential appearing in the integration of the definition, and the function is expected to be quite badly behaved as the examples show.

The graph of $S_b(x)$ for $b=0.7$ is given in Figure \ref{Sx07}:
\begin{figure}[H]
\centering
\includegraphics[width=120mm]{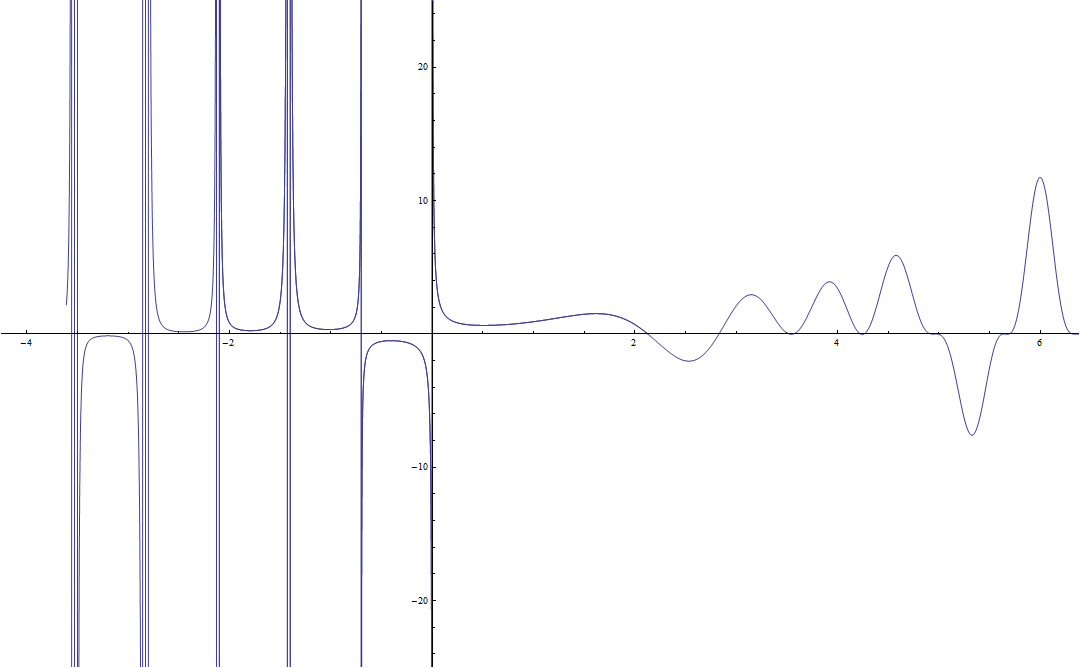}
\caption{The graph of $S_b(x)$ for $b=0.7, -2Q<x<3Q$, Q=2.129}
\label{Sx07}
\end{figure}

The graph of $S_b(x)$ for $b=0.248$ is given in Figure \ref{Sx0248}:
\begin{figure}[H]
\centering
\includegraphics[width=120mm]{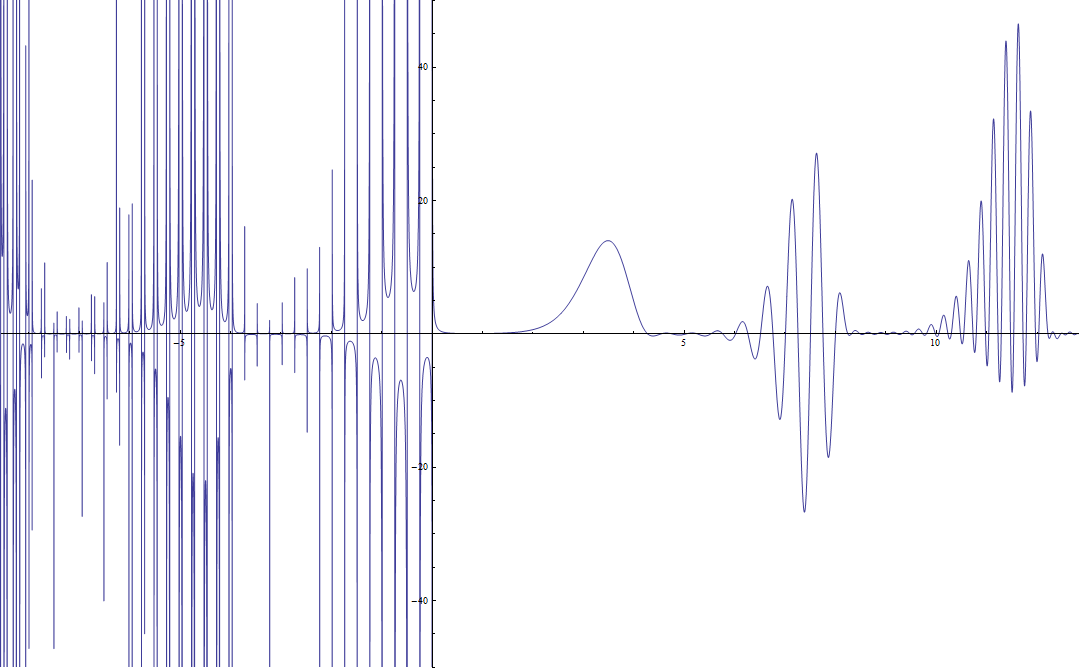}
\caption{The graph of $S_b(x)$ for $b=0.248, -2Q<x<3Q$, Q=4.280}
\label{Sx0248}
\end{figure}

The graph of $S_b(x)$ for $b=0.1$ is given in Figure \ref{Sx01}:
\begin{figure}[H]
\centering
\includegraphics[width=120mm]{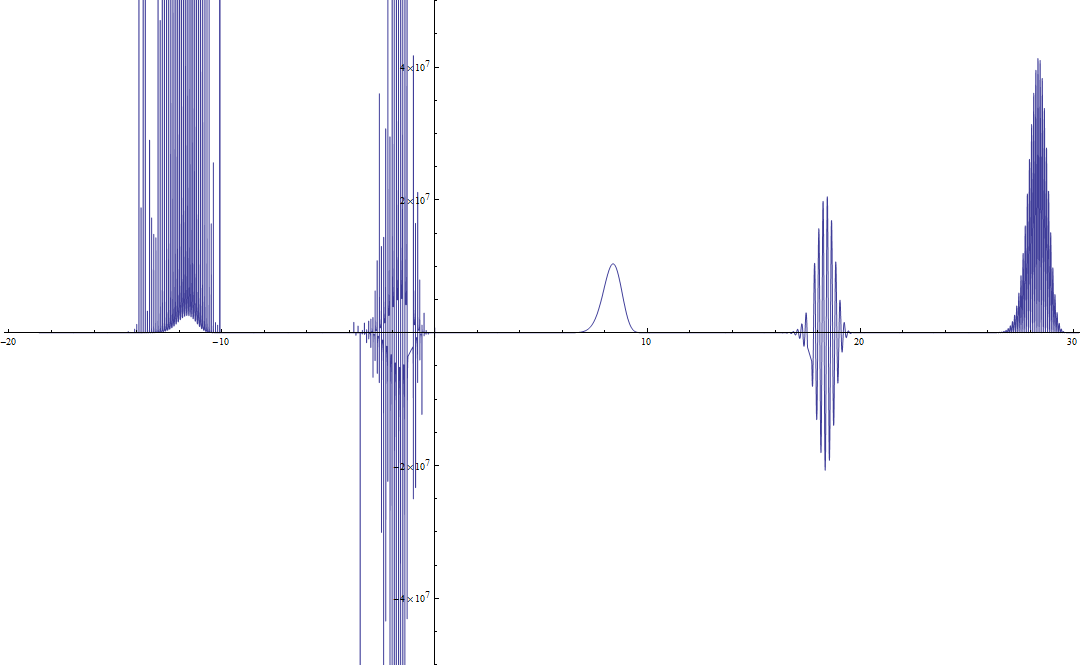}
\caption{The graph of $S_b(x)$ for $b=0.1, -2Q<x<3Q$, Q=10.1}
\label{Sx01}
\end{figure}

\begin{figure}[H]
\centering
\includegraphics[width=140mm]{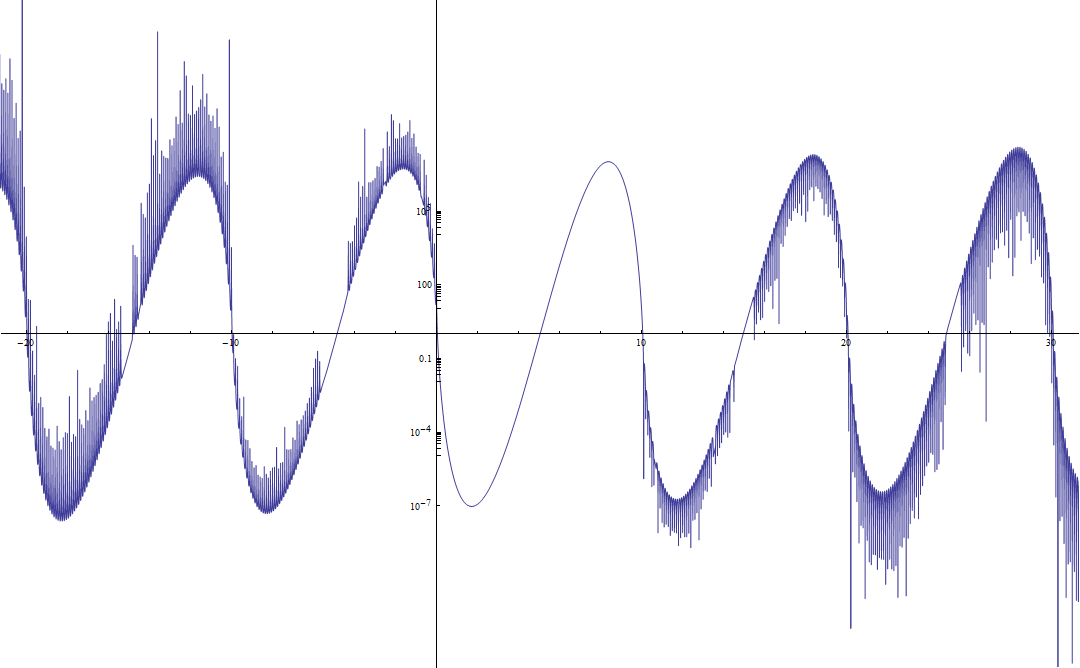}
\caption{The logarithmic plot of $|S_b(x)|$ for $-2Q<x<3Q$ and $b=0.1$}
\label{logSx}
\end{figure}

The case for $b=0.1$ can be better illustrated by taking the log plot (of its absolute values). From Figure \ref{logSx}, the spikes correspond to the position of the poles and zeroes.

Again let us comment on certain features of the graphs
\begin{itemize}
\item We can clearly observe the locations of the poles and zeroes. The zeroes are simple from the fact that the values of the function should alternate in signs. For $b=0.7$, the appearance of a seemingly critical points comes from the fact that $b\inv = 1.429$ is quite close to $2b=1.4$ as explained before.
\item By the functional equation \eqref{funceq}, the whole graph really depends on the portion $0<x<Q$. The size of the bump governs the recurring amplitude of the graph by successive factor of 2, which comes from $2\sin(\pi b)\sin(\pi b\inv)$ appearing in the reflection properties for $S_b(x)$. This also explains the "periodicity" of the log plot for $b=0.1$, which, in fact, has shifted by $\ln(2)$.
\item The portion of the graph $0<x<Q$ has exactly one local minimum and one local maximum for all values of $b$, and they are all positive real numbers. By taking the log plot in Figure \ref{logSQx}, which essentially is the value of the integral, (we rescaled so that the endpoints at $x=0$ and $x=Q$ matches), we see immediately that the size of the bump has a symmetric exponential relation along $x=\frac{Q}{2}$.
 \begin{figure}[H]
\centering
\includegraphics[width=140mm]{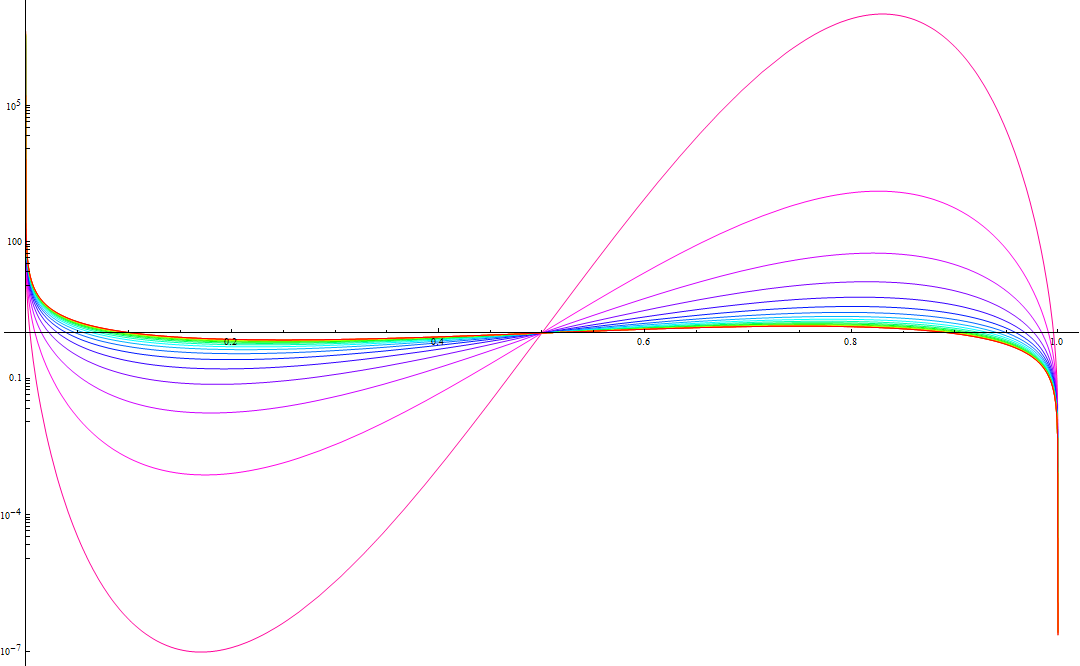}
\caption{The logarithmic plot of $S_b(Qx)$ for $0<x<1$ and $b=0.1$ to $1$ in steps of $0.05$}
\label{logSQx}
\end{figure}
\end{itemize}

\subsection{Values along imaginary direction}
The values along the imaginary direction are simply governed by the asymptotic behavior. For completeness, we illustrate the graphs for $G_b(ix),G_b(\frac{Q}{2}+ix)$ and $G_b(Q+ix)$ which appears quite often in the calculation of the representation theory of the quantum plane \cite{Ip2}. Again we choose $b=0.7$. The absolute value is represented by a thick line, real value by a blue line, and imaginary value by a red line.

\begin{figure}[H]
\centering
\includegraphics[width=120mm]{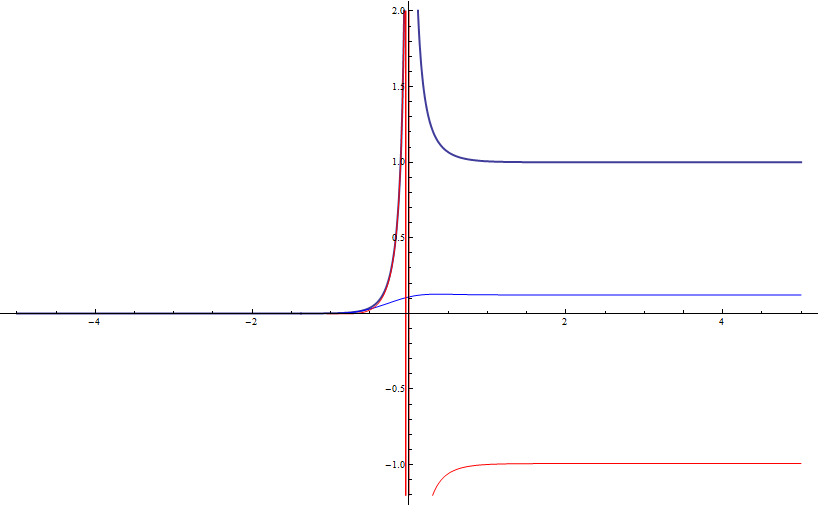}
\caption{The graph of $G_b(ix)$ at $b=0.7$ for $-5<x<5$}
\end{figure}

\begin{figure}[H]
\centering
\includegraphics[width=120mm]{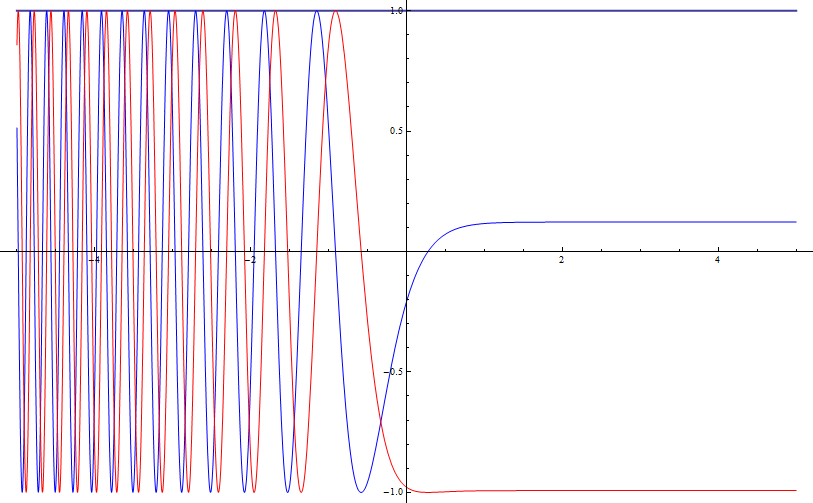}
\caption{The graph of $G_b(\frac{Q}{2}+ix)$ at $b=0.7$ for $-5<x<5$}
\end{figure}

\begin{figure}[H]
\centering
\includegraphics[width=120mm]{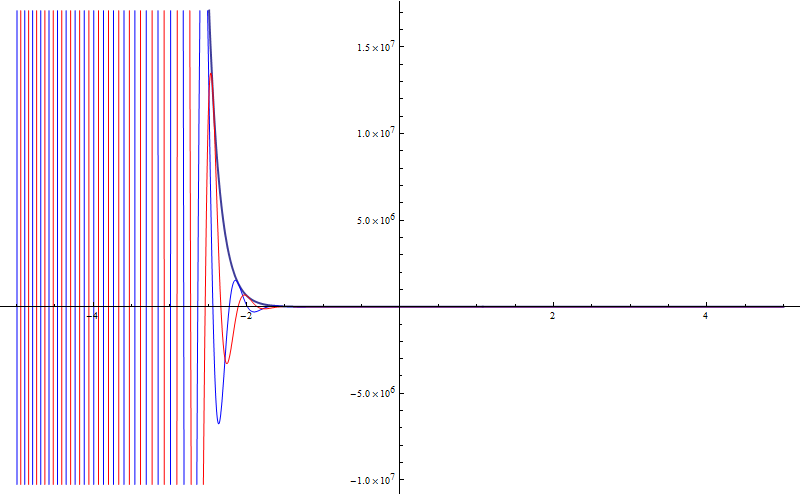}
\caption{The graph of $G_b(Q+ix)$ at $b=0.7$ for $-5<x<5$}
\end{figure}

\begin{itemize}
\item We observe the asymptotic behavior, that as $x\to \oo$, they all approach the complex number $\over[\ze_b]=0.123-0.992i$.
\item For $x\to-\oo$, $|G_b(ix)|\to e^{\pi Qx}$, while $|G_b(Q+ix)|\to e^{-\pi Qx}$. The relationship of the two graphs can easily be derived from the reflection properties \eqref{reflection}.
\item  Finally $|G_b(\frac{Q}{2}+ix)|=1$, and the asymptotic behavior for $x\to-\oo$ is $\overline{\ze_b}e^{-\pi i\frac{Q^2}{4}}e^{-\pi ix^2}$, hence the increase in frequency is governed by $\pi x^2$.
\end{itemize}

\section{The Compact Quantum Dilogarithm}\label{sec:compact}
\subsection{Visualizations}
The analytic properties extend to the compact case, when $Im(b^2)>0$, using an infinite product expansion, which we will call it the \emph{compact quantum dilogarithm}:

\begin{Prop}\label{infprod}\cite[Prop 5]{Sh} For $Im(b^2)>0$, we have the infinite product descripiton
\Eq{G_b(z)=\bar{\ze_b}\frac{\prod_{n=1}^\oo (1-e^{2\pi ib\inv(x-nb\inv)})}{\prod_{n=0}^\oo(1-e^{2\pi ib(x+nb)})}.}
\end{Prop}

The location of the poles and zeroes are still the same, with poles at $-nb-mb\inv$, and zeroes at $Q+nb+mb\inv$ for $n,m=0,1,2,...$

\begin{figure}[H]
\centering
\includegraphics[width=120mm]{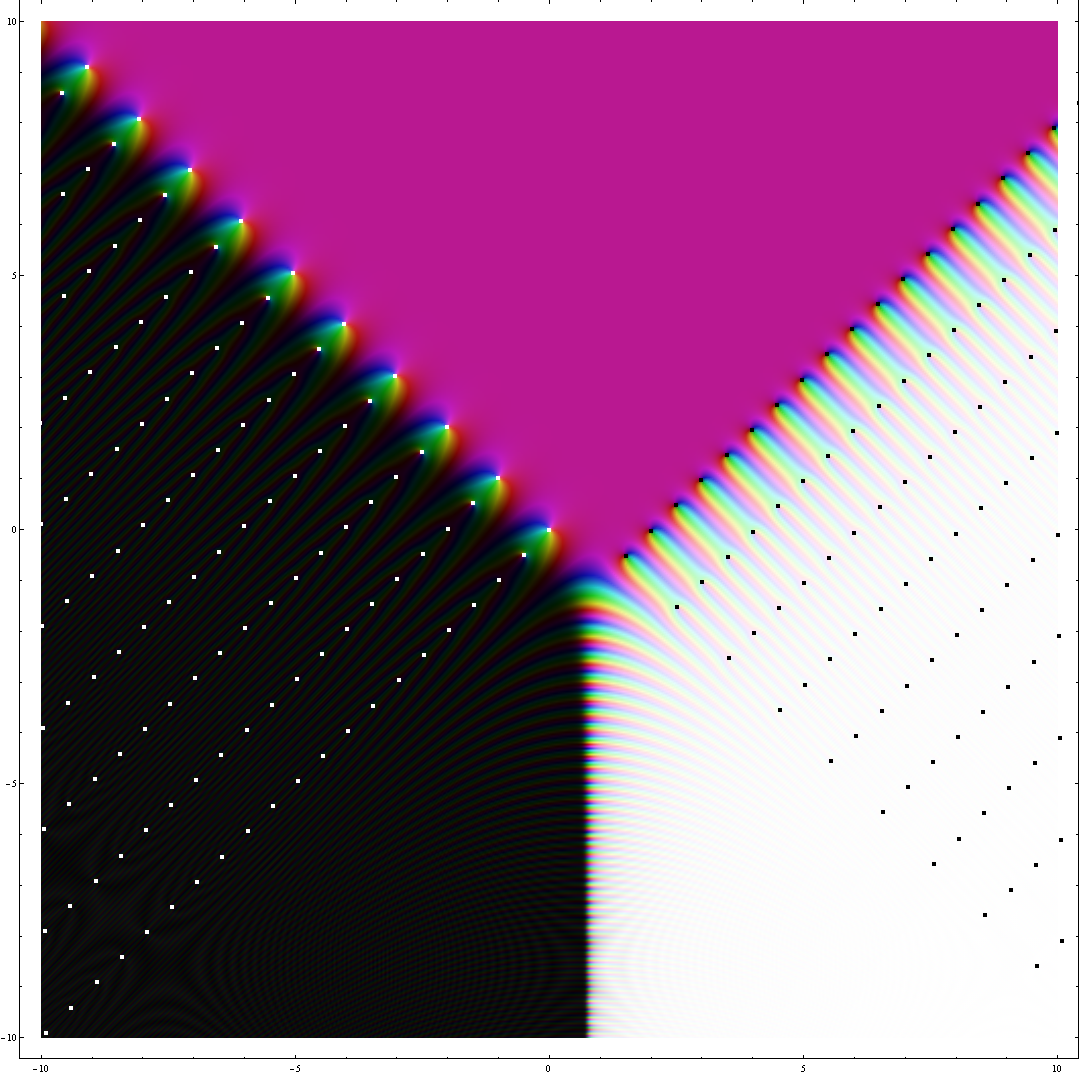}
\caption{The graph of $G_b(z)$ at $b=e^{\pi i/4}0.7$ for $z\in[-10,10]\x[-10,10]$}
\label{Gx07i}
\end{figure}

The most interesting case occurs when $b^2=ir$ for $r>0$, in which $0<q=e^{\pi ib^2}<1$. We establish the graph for $b=e^{\pi i/4}0.7$ in Figure \ref{Gx07i}

We have highlighted the zeroes with black dots and poles with white dots that align nicely together. We observe that the exponential decay and growth are extremely fast, due to the asymmetry bringing in a real part for the asymptotic $e^{\pi iz(z-Q)}$, thus giving a quadratic exponential growth.

\begin{figure}[H]
\centering
\includegraphics[width=120mm]{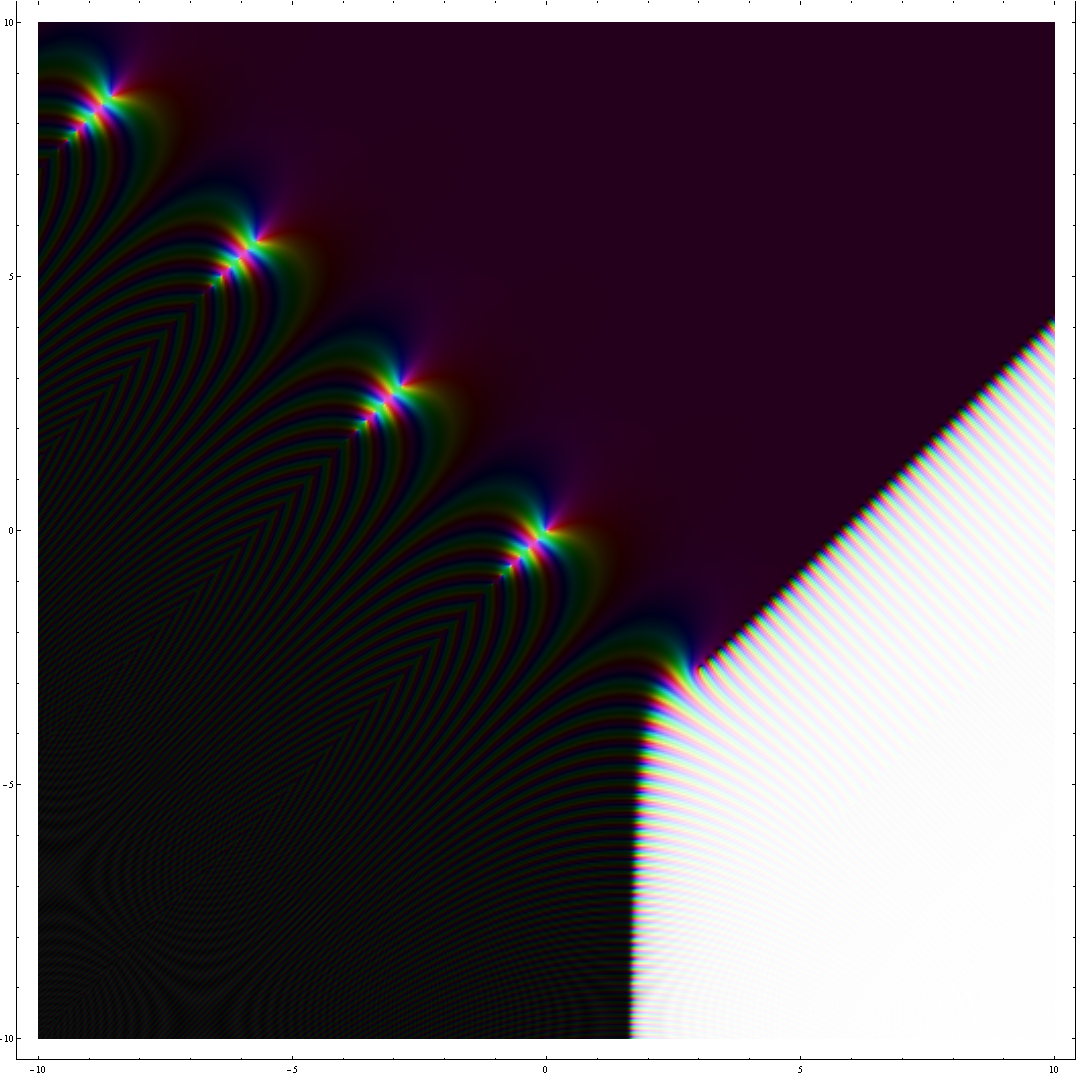}
\caption{The graph of $G_b(z)$ at $b=e^{\pi i/4}0.248$ for $z\in[-10,10]\x[-10,10]$}
\label{Gx0248i}
\end{figure}

The case for $b=e^{\pi i/4}0.248$ illustrated in Figure \ref{Gx0248i} is less interesting due to the quick growth and decay of the function.

We observe that the separation of the poles is stretching out with increased exponential decay/growth compared with bigger $|b|$.

\subsection{Limit of $\til[G_b](z)$ to $\G(z)$ as $b\to0$}

Finally let us look at the theorem proved in \cite{Ip1}:
\begin{Thm}\label{limit} The following limit holds for $b^2=ir\to i0^+$
\Eq{\label{Glim}\lim_{r\to 0} \frac{G_b(bz)}{\sqrt{-i}|b|  (1-q^2)^{x-1}} = \G(z),}
where $\sqrt{-i}=e^{-\frac{\pi i}{4}}$ and $-\frac{\pi}{2}<\arg(1-q^2)<\frac{\pi}{2}$.
\end{Thm}

Let us demonstrate the limiting process by considering the graph of the above function 
\Eq{\til[G]_b(z):=\frac{G_b(bz)}{\sqrt{-i}|b|  (1-q^2)^{z-1}}.}

\begin{figure}[H]
\centering
\includegraphics[width=120mm]{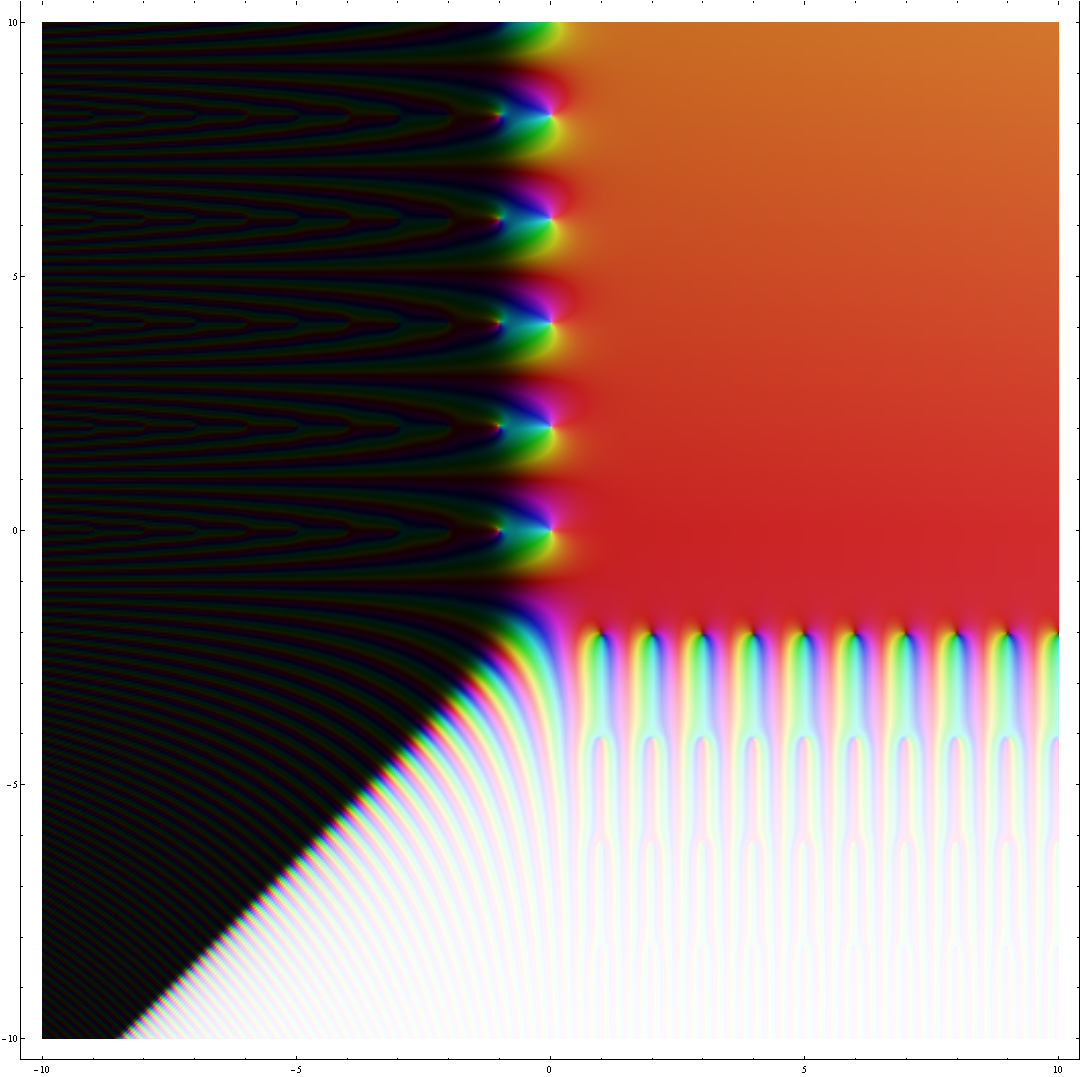}
\caption{The graph of $\til[G]_b(z)$ at $b=e^{\pi i/4}0.7$ for $z\in[-10,10]\x[-10,10]$}
\label{GG07}
\end{figure}

From Figure \ref{GG07}, we see that under multiplication by $b$ in the argument, the zeroes and poles are aligned with the axis. The first zero occur at $z=\frac{Q}{b}=1-\frac{1}{|b|^2}i$. Now as $|b|\to 0$, we see that $\frac{Q}{b}$ moves downward and pushes away the zeroes and poles. It becomes clear in the graph for smaller $|b|$.

\begin{figure}[H]
\centering
\includegraphics[width=120mm]{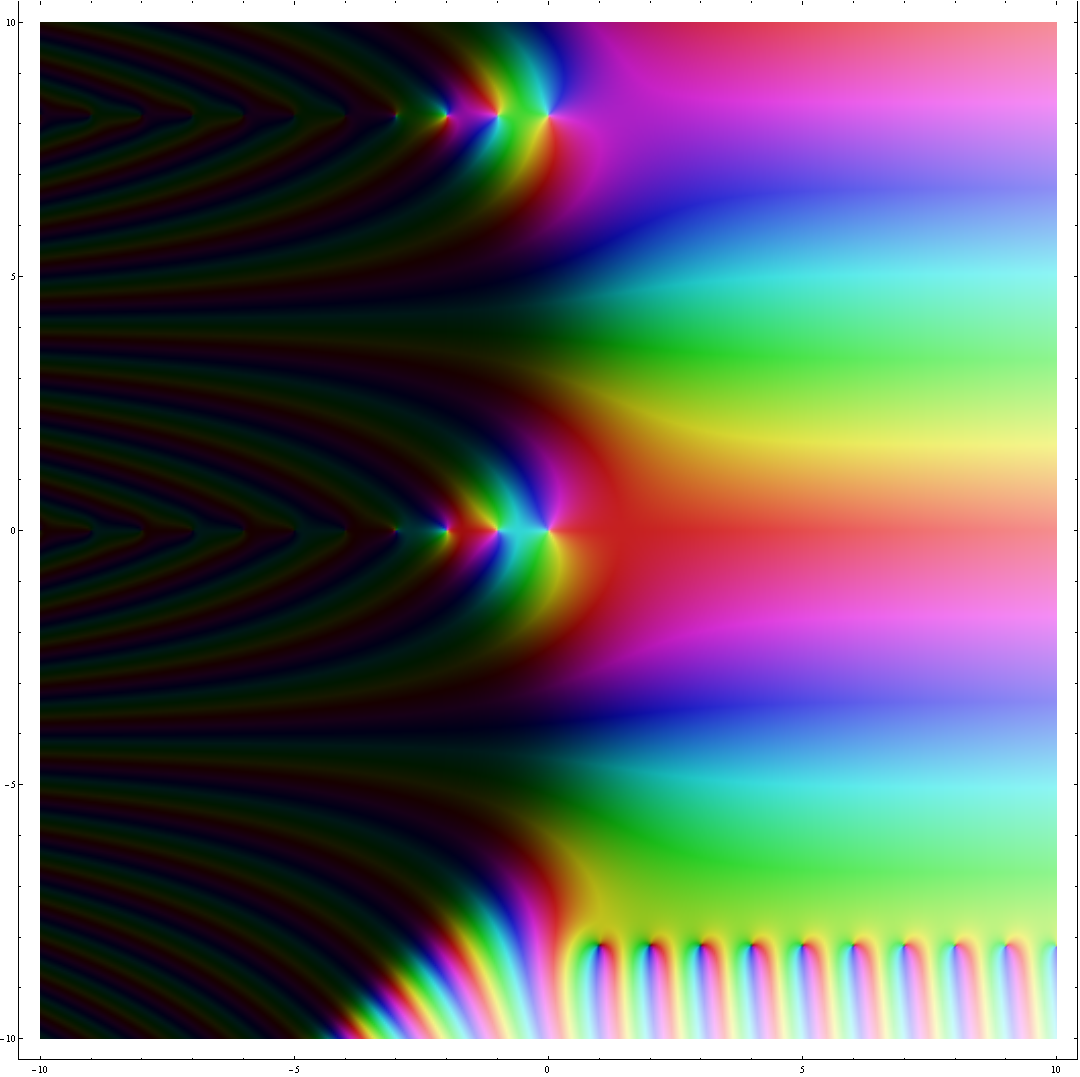}
\caption{The graph of $\til[G]_b(z)$ at $b=e^{\pi i/4}0.35$ for $z\in[-10,10]\x[-10,10]$}
\label{GG035}
\end{figure}

\begin{figure}[H]
\centering
\includegraphics[width=120mm]{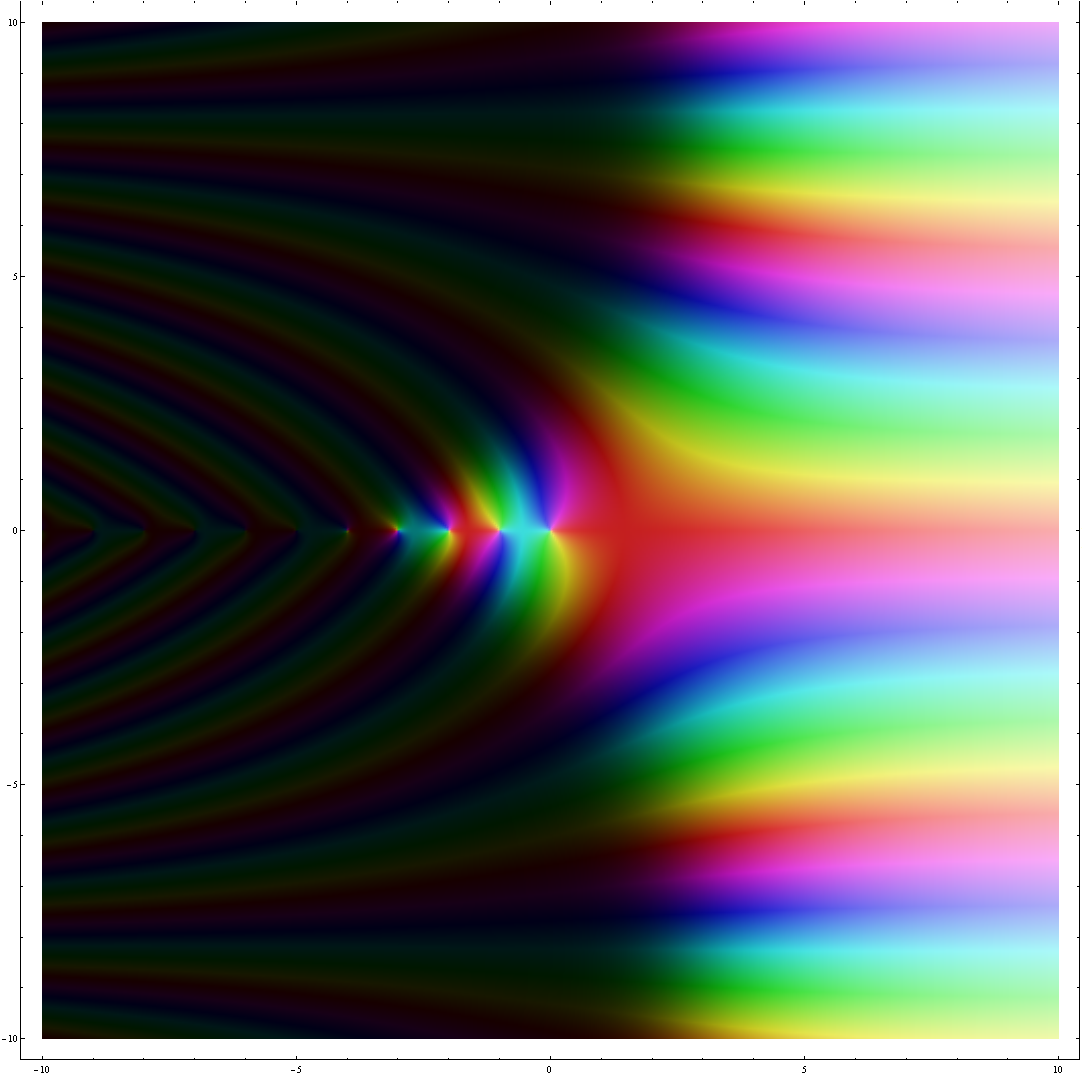}
\caption{The graph of $\til[G]_b(z)$ at $b=e^{\pi i/4}0.248$ for $z\in[-10,10]\x[-10,10]$}
\label{GG0248}
\end{figure}

Now as $|b|$ becomes smaller, from Figure \ref{GG035} we see that the small area around center begins to take shape with the poles and zeroes pushed away. Further reducing $|b|$ so that $\frac{Q}{b}$ is out of plotting range, we can see immediately in Figure \ref{GG0248} that indeed the graph resembles the Gamma function $\G(z)$ given in Figure \ref{Gamma}.

\section{Relating $G_b(z)$ to other variants}\label{sec:relation}
We will see below that all variants of $G_b$ amounts to shifting, scaling, taking inverses and multiplication by a constant. Hence the graphs of these variants can be obtained from those of $G_b(z)$ easily by shifting, stretching, and color phase changes.

The important function that satisfies the pentagon equation for Weyl-type operator is defined by
\Eq{g_b(z)=\frac{\over[\ze_b]}{G_b(\frac{Q}{2}+\frac{\log z}{2\pi ib})}.}
\begin{Lem} Let $u, v$ be positive self-adjoint operators with $uv=q^2vu$, $q=e^{\pi i b^2}$. Then
\Eq{\label{qexp}g_b(u)g_b(v)=g_b(u+v),}
\Eq{\label{qpenta}g_b(v)g_b(u)=g_b(u)g_b(q\inv uv)g_b(v).}
\end{Lem}

However, the graph of $g_b$ which carries a branch cut at $x\in(-\oo,0)$ is not particularly interesting, due to slow growth, and the fact that $-\pi < Im(\log z)<\pi$, so that 
\Eq{0<\frac{Q}{2}-\frac{1}{2b}<Re(\frac{Q}{2}+\frac{\log z}{2\pi ib}) < \frac{Q}{2}+\frac{1}{2b}<Q,}
and hence the argument in $G_b$ never hits any zeroes or poles.

\begin{figure}[H]
\centering
\includegraphics[width=60mm]{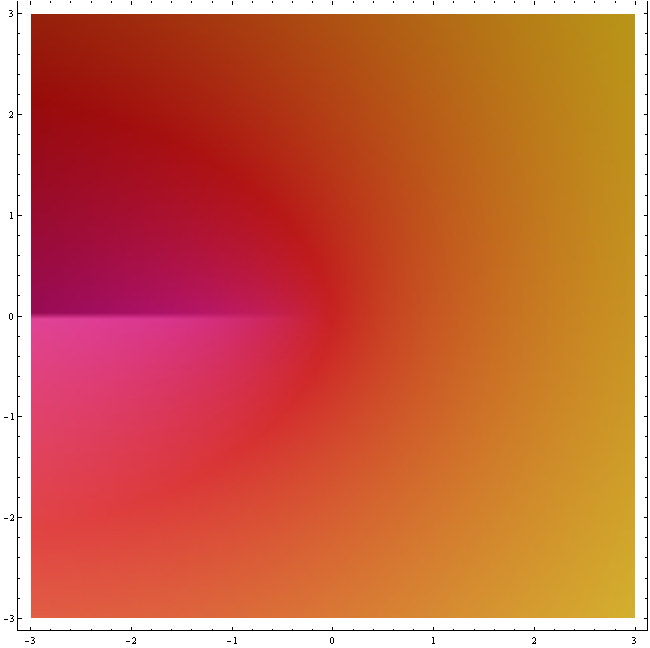}
\caption{The graph of $g_b(z)$ for $b=0.7$ for $z\in[-3,3]\x[-3,3]$}
\label{gb}
\end{figure}

The quantum dilogarithm $G_ b(x)$ defined above has real $b$ as parameters. Although the integral formula can be analytically continued to complex $b$ with positive real part, we recall that the function coincide with the \emph{compact quantum dilogarithm} given in Prop \ref{infprod}:
\Eq{G_b(z)=\bar{\ze_b}\frac{\prod_{n=1}^\oo (1-e^{2\pi ib\inv(x-nb\inv)})}{\prod_{n=0}^\oo(1-e^{2\pi ib(x+nb)})}.}

From this we can compare with $\c(z)$ defined by Volkov \cite{Vo}:
\begin{Def} Let $q=e^{\pi i\t}$ for $Im(\t)>0$. Then the hyperbolic Gamma function is defined by
\Eq{\c(z):=\frac{(q^2e^{-2\pi iz};q^2)_\oo}{(e^{-2\pi iz/\t};q^{-2/\t^2})_\oo}.}
\end{Def}
\begin{Prop} With $\t=b^2$, we have 
\Eq{\c(z)=\over[\ze_b]G_b(Q-\frac{z}{b})\inv.}
\end{Prop}

The infinite product expansion is also involved in the original definition $\psi_b$ by Faddeev \cite{FKa}:
\begin{Def} Let $s_q(w)=\prod_{n=0}^\oo (1+q^{2n+1}w)$ and define 
\Eq{\psi(p):=\frac{s_q(e^{bp})}{s_{\til[q]}(e^{p/b})},} where $\til[q]=e^{\pi ib^{-2}}$. Then $\psi(p)$ admits the integral expansion
\Eq{\psi(p):=\exp\left(\frac{1}{4}\int_\W \frac{e^{ip\xi/\pi}}{\sinh(b\xi)\sinh(\xi/b)}\frac{d\xi}{\xi}\right).}
\end{Def}

Fock-Goncharov \cite{FG} used a similar variant $\Phi^\hbar$, which is the same as $V_\h$ used by Woronowicz-Zakrzewski's \cite{WZ}:
\begin{Def}
\Eq{\Phi^{\hbar}(z):=\exp\left(-\frac{1}{4}\int_\W \frac{e^{-ipz}}{\sinh(\pi p)\sinh(\pi \hbar p)}\frac{dp}{p}\right),}
\Eq{V_\h(z):=\exp\left(\frac{1}{2\pi i}\int_0^\oo \log(1+x^{-\h})\frac{dx}{x+e^{-z}}\right).}
\end{Def}
The we have the relations
\begin{Prop}
\Eq{V_{b^{-2}}(z)=\Phi^{b^2}(z)=g_b(e^{z})\inv=\ze_bG_b\left(\frac{Q}{2}-\frac{iz}{2\pi b}\right),}
\Eq{\psi(z)=\Phi^{b^2}(-bz)\inv=\c(b(Q+\frac{z}{2\pi i}))=\over[\ze_b]G_b(\frac{Q}{2}+\frac{iz}{2\pi})\inv.}
\end{Prop}

Finally, a more general version with two parameters is given by Ruijinaars \cite{Ru}
\begin{Lem} Let $a_-,a_+$ be two positive real numbers. Ruijinaars's G function is defined by
\Eq{G(a_+,a_-;z):=\exp\left(i\int_0^\oo \frac{dy}{y}\left(\frac{\sin 2yz}{2\sinh(a_+y)\sinh(a_-y)}-\frac{z}{a_+a_-y}\right)\right),}
with $|Im(z)|<(a_++a_-)/2$ and extend meromorphically to the whole complex plane.

Then \Eq{G(b,b\inv;z)=S_b(\frac{Q}{2}-iz).}
\end{Lem}

\end{document}